\documentclass[12pt]{article}
\pdfoutput=1
\usepackage{graphicx,pict2e,latexsym,xcolor}
\usepackage[square,numbers]{natbib}
\newtheorem{defn}{Definition}
\newtheorem{theorem}{Theorem}
\newtheorem{lemma}{Lemma}

\def\dot{\small $\bullet$}
\def\rdot{\color{red}\small $\bullet$}

\def\gdot{\color{green}\small $\bullet$}

\begin{document}
\bibliographystyle{plainnat}
\pagestyle{plain}

\title{\Large \bf Approximate Information Tests \\ on Statistical Submanifolds}

\author{Michael W. Trosset\thanks{Department of Statistics, Indiana University.  
E-mail: {\tt mtrosset@indiana.edu}} \and
Carey E. Priebe\thanks{Department of Applied Mathematics \& Statistics, Johns Hopkins University. 
E-mail: {\tt cep@jhu.edu}} \and
}

\date{\today}

\maketitle


\begin{abstract}
Parametric inference posits a statistical model that is a specified family of probability distributions.
Restricted inference, e.g., restricted likelihood ratio testing,
attempts to exploit the structure of a statistical submodel that is a subset of the specified family.  We consider the problem of testing a simple hypothesis against alternatives from such a submodel.  In the case of an unknown submodel, it is not clear how to realize the benefits of restricted inference.
To do so, we first
construct information tests that are locally asymptotically equivalent to likelihood ratio tests.  Information tests are conceptually appealing but (in general) computationally intractable.  However, unlike restricted likelihood ratio tests, restricted information tests can be approximated even when the statistical submodel is unknown.
We construct approximate information tests using manifold learning procedures to extract information from samples of an unknown (or intractable) submodel,
thereby providing a roadmap for computational solutions to a class of previously impenetrable problems in statistical inference.  Examples illustrate the efficacy of the proposed methodology.
\end{abstract}


\newpage

\tableofcontents

\newpage


\section{Introduction}
\label{intro}

An engrossing challenge arises when an appropriate statistical model is a subset of a familiar family of probability distributions:
how to exploit the structure of the restricted model for the purpose of subsequent inference?
This challenge encompasses theoretical, methodological, computational, and practical concerns.
The reasons to address these concerns
are especially compelling when the restricted model is of lower dimension than the unrestricted model,
as parsimony principles encourage the selection of less complicated models.

The following example illustrates the concerns of the present manuscript.

\subparagraph{Motivating Example}
Consider a multinomial experiment with $7$ possible outcomes and probability vector 
$\theta \in \Re^7$.  To test the simple null hypothesis
\[
H_0: \theta = \bar{\theta} = (0.09,0.09,0.09,0.25,0.16,0.16,0.16)
\]
at significance level $\alpha=0.05$, we perform $n=30$ trials and observe
\[
o = (3,5,4,6,9,2,1).
\]
Should we reject $H_0$?

The likelihood ratio test statistic of
\[
G^2 = 2 \sum_{j=1}^7 o_j \log \left( o_j/n\bar{\theta}_j \right) = 11.93649
\]
results in an (approximate) significance probability of $\mbox{\bf p} = 0.0634$.  Pearson's $X^2 = 11.23519$ results in $\mbox{\bf p} = 0.0814$.  Neither test provides compelling evidence against $H_0$.

Suppose, however, that it is possible to perform an auxiliary experiment that randomly generates possible values of $\theta$ for the primary experiment.  The auxiliary experiment is performed $m=100$ times and it is found that $96$\% of the variation in the $m=100$ values of $\theta$ is explained by $2$ principal components.  This finding suggests the possibility that $\theta$ is restricted to a (slightly curved) $2$-dimensional submanifold of the $6$-dimensional simplex.  Can this revelation be exploited to construct a more powerful test?

If the submanifold was known, then one could perform a restricted likelihood ratio test.  But {\em the submanifold is not known}.  \hfill $\Box$

\bigskip

In fact, the family of multinomial distributions provides numerous examples of dimension-restricted submodels.
In statistical genetics, the phenomenon of
Hardy-Weinberg equilibrium corresponds to a much-studied $1$-parameter subfamily of trinomial distributions.
Spherical subfamilies of multinomial distributions
\citep{gous:1999}
are potentially valuable in a variety of applications, e.g., text mining \citep{hall&hoffman:2000}.
In a recent effort to discover brainwide neural-behavioral maps from optogenetic experiments on {\em Drosophila}\/ larvae
\citep{vogelstein&etal:2014},
each neuron line was modeled by a $29$-dimensional vector of multinomial probabilities but
the available evidence suggested that these vectors resided on an unknown $4$-dimensional submanifold.
These examples suggest a natural progression,
from a submodel that is known and tractable,
to a submodel that is known but possibly intractable,
to an unknown submodel that can be sampled,
to an unknown submodel that must be estimated.
The particular challenge of how to exploit low-dimensional structure that is apparent but unknown motivated our investigation.  The present manuscript addresses the case of known submodels and unknown submodels that can be sampled; a sequel will address the case of unknown submodels that must be estimated.

\begin{figure}[tb]
\begin{center}
\fbox{ 
\begin{minipage}{5.5in}
\vspace{1em}
Suppose that known distributions $\bar{p},p_1,\ldots,p_m$ lie in an unknown statistical submanifold.
To test $H_0: p = \bar{p}$ against alternatives that lie in the submanifold, we propose the following procedure.
\begin{enumerate}

\item Compute $h_{ij}$, the pairwise Hellinger distances between 
$\bar{p},p_1,\ldots,p_m$.

\item Construct ${\mathcal G}$, a graph whose vertices correspond to the known distributions.  Connect vertices $i$ and $j$ when $h_{ij}$ is sufficiently small.

\item Compute the pairwise shortest path distances in 
${\mathcal G}$.

\item Construct $\bar{z},z_1,\ldots,z_m \in \Re^r$, an embedding of 
${\mathcal G}$ whose
pairwise Euclidean distances approximate the pairwise shortest path distances.

\item From $x_1,\ldots,x_n \sim p$, construct a nonparametric density estimate $\hat{p}_n$.  Compute the Hellinger distances of 
$\hat{p}_n$ from $p_1,\ldots,p_m$ and 
embed $\hat{p}_n$ as $y(\vec{x}) \in \Re^r$
in the previously constructed Euclidean representation.  The proposed test rejects $H_0: \theta=\bar{\theta}$ if and only if the test statistic
$\left\| y \left( \vec{x} \right) -\bar{z} \right\|$
is sufficiently large.

\item Estimate a significance probability by generating simulated random samples from the hypothesized distribution $\bar{p}$.

\end{enumerate}
\vspace{1em}
\end{minipage}
}
\end{center}
\caption{An approximate information test for the case of an unknown submodel that can be sampled.  Steps 2--4 are essentially Isomap \citep{isomap:2000}, used here to represent the Riemannian structure of a statistical manifold rather than a data manifold.  Details are provided in Section~\ref{AIT}.}
\label{fig:AIT}
\end{figure}

For unknown submodels that can be sampled, we propose the computationally intensive {\em approximate information test}\/ summarized in Figure~\ref{fig:AIT}.  The theory that underlies and motivates this procedure originates in information geometry, specifically in the well-known fact that Fisher information induces Riemannian structure on a statistical manifold.
It leads to {\em information tests}\/ that are conceptually appealing but (in general) computationally intractable.  Approximate information tests circumvent the intractability of information tests.

Sections~\ref{prelim}--\ref{HW} develop and illustrate the theory of information tests.
Section~\ref{prelim} establishes the mathematical framework that informs our investigation.  We review the fundamental concepts of a statistical manifold and the Riemannian structure induced on it by Fisher information.  We demonstrate that information distance, i.e., geodesic distance on this Riemannian manifold, is more practically derived from Hellinger distance, and we briefly review minimum Hellinger distance estimation.
Sections~\ref{IT}--\ref{HW} develop tests of simple null hypotheses using the concept of information distance.
Section~\ref{IT} demonstrates that information tests are locally asymptotically equivalent to various classical tests (Hellinger distance, Wald, likelihood ratio, and Hellinger disparity distance).  Section~\ref{RIT} derives information tests for submodels of the multinomial model.
Section~\ref{HW} provides examples using the Hardy-Weinberg submodel of the trinomial model.

Despite their conceptual appeal, the information tests developed in Sections~\ref{IT}--\ref{HW} are of limited practical application.
Hence, our primary contribution lies in Section~\ref{AIT}, which proposes a discrete approximation of an information test and illustrates its effectiveness in two cases for which an unknown submodel can be sampled.  Section~\ref{disc} discusses implications and possible extensions.

\section{Preliminaries} 
\label{prelim}

\subsection{Statistical Manifolds}
\label{manifold}

We begin by recalling some basic properties of differentiable manifolds.
See \citep{matsushima:1972} for a more detailed explication of these concepts.
Let $M$ denote a completely separable Hausdorff space.  
Let $U \subseteq M$ and $V \subseteq \Re^k$ denote open sets.
If $\varphi : U \rightarrow V$ is a homeomorphism,
then $\varphi(u) = (x_1(u),\ldots,x_k(u))$ defines a coordinate system on $U$.
The $x_i$ are the coordinate functions and $\varphi^{-1}$ is a parametrization of $U$.  The pair $(U,\varphi)$ is a chart.
An atlas on $M$ is a collection of charts $\{ (U_a,\varphi_a) \}$ such that the $U_a$ cover $M$.

The set $M$ is a $k$-dimensional topological manifold if and only if it admits an atlas for which each $\varphi_a(U_a)$ is open in $\Re^k$.
It is a differentiable manifold if and only if the transition maps $\varphi_b \varphi_a^{-1}$ are diffeomorphisms.
A subset $S \subset M$ is a $d$-dimensional embedded submanifold if and only if, for every $p \in S$, there is a chart $(U,\varphi)$ such that $p \in U$
and
\[
\varphi(U \cap S) = \varphi(U) \cap \left( \Re^d \times 
\{ \vec{0} \in \Re^{k-d} \} \right) =
\left\{ y \in \varphi(U) : y_{d+1}= \cdots = y_k=0 \right\}.
\]

Our explication of statistical manifolds follows  \citet{murray&rice:1993}, from whom much of our notation is borrowed.  
Let $(\Omega,\mathcal{B},\mu)$ denote a measure space.
Let $\mathcal{M}$ denote the nonnegative measures on $(\Omega,\mathcal{B})$ that are absolutely continuous with respect to $\mu$.
We write an element of $\mathcal{M}$ as $p \, d\mu$, where $p$ is a density function with respect to $\mu$.
We write $p \, d\mu \sim q \, d\mu$ and say that $p \, d\mu$ and $q \, d\mu$ are equivalent up to scale if and only if 
\[
\frac{\int_B p(x) \, d\mu(x)}{\int_\Omega p(x) \, d\mu(x)} =
\frac{\int_B q(x) \, d\mu(x)}{\int_\Omega q(x) \, d\mu(x)}
\]
for every $B \in \mathcal{B}$.
\citet{murray&rice:1993} regard a probability measure as an equivalence class of finite measures.  
Let $\mathcal{P}$ denote the space of probability measures in $\mathcal{M}$, i.e., the set of finite measures up to scale.

Let $\Re_\Omega$ denote the vector space of measurable real-valued functions on $\Omega$ and define the log-likelihood map
$\ell : \mathcal{M} \rightarrow \Re_\Omega$ by
$\ell(p \, d\mu) = \log(p)$.  We say that the log-likelihood map is smooth
if and only if, for each $x \in \Omega$, the corresponding real-valued component map defined by $p \, d\mu \mapsto [\log(p)](x)$ is sufficiently differentiable.

\begin{defn}
\label{df:statmanifold}
Let 
$P = \left\{ p(\cdot,\theta) \, d\mu \; : \; \theta \in \Theta \subseteq \Re^k \right\}$
denote a parametric family of probability distributions in $\mathcal{P}$.
We say that $P$ is a statistical manifold if and only if $P$ is a differentiable manifold, the log-likelihood map is smooth, and, for any $p \, d\mu \in P$, the random variables
\[
\frac{\partial \ell}{\partial \theta^1} \left( p \, d\mu \right),
\ldots,
\frac{\partial \ell}{\partial \theta^k} \left( p \, d\mu \right)
\]
are linearly independent.
\end{defn}

We might dispense with the parametric structure of $P$, but many of the familiar concepts and results of classical statistics are stated in terms of index sets rather than families of distributions.  For example,
fix $p \, d\mu \in P$.  Then the random vector
\[
d_p\ell = \left(
\frac{\partial \ell}{\partial \theta^1} \left( p \, d\mu \right),
\ldots,
\frac{\partial \ell}{\partial \theta^r} \left( p \, d\mu \right)
\right)
\]
is the score vector at $p \, d\mu$,
and the set of vectors obtained by observing the score vector at each $x \in \Omega$ is the tangent space of $P$ at $p \, d\mu$, denoted $T_pP$.
Our exposition will emphasize the manifold structure of $P$ itself, but one can just as easily regard $P$
as indexed by a $k$-dimensional manifold $\Theta$---and it is often convenient to do so.

\subsection{Riemannian Geometry and Fisher Information}
\label{I}

A metric tensor on the statistical manifold $P$ is a collection of inner products on the tangent spaces of $P$.  If $P$ admits a metric tensor, then $P$ is a Riemannian manifold.  See \citep[Part II]{milnor:1963} and \citep{hicks:1971} for concise introductions to Riemannian geometry.  Note that many authors refer to the metric tensor as a Riemannian metric.  In neither case is the word ``metric'' used in the sense of a distance function.

Let $E_p$ denote expectation with respect to $p \, d\mu$, i.e.,
$E_p f = \int_\Omega f(x) p(x) \, d\mu(x)$.
Define an inner product on the space of square-integrable $f$ by
$\langle f,g \rangle_p = E_p fg$.
If the log-likelihood map is smooth,
then the Fisher information matrix $I(p) = [g_{ij}(p)]$ has entries
\[
g_{ij}(p) = E_p \frac{\partial \ell}{\partial \theta^i}
\frac{\partial \ell}{\partial \theta^j} =
\left\langle \frac{\partial \ell}{\partial \theta^i},
\frac{\partial \ell}{\partial \theta^j} \right\rangle_p.
\]
Because the scores are linearly independent, $I(p)$ is the matrix of the inner product $\langle \cdot,\cdot \rangle_p$ with respect to the basis defined by the scores.
 
\citet{rao:1945} observed that Fisher information induces a natural metric tensor on $P$.  To obtain a coordinate-free representation of this tensor, i.e., a representation that does not involve Fisher information matrices, suppose that $v \in T_pP$ and let $\gamma : (-\epsilon,\epsilon) \rightarrow P$ be any variation with tangent vector $v$ at $p=\gamma(0)$.  The differential of the log-likelihood map at $p$ is the function $d_p\ell : T_pP \rightarrow \Re_\Omega$ defined by
\[
d_p\ell(v) = (\ell \circ \gamma)^\prime(0) =
\frac{d}{dt} \ell(\gamma(t))|_{t=0} =
\lim_{t \rightarrow 0} \frac{\ell(\gamma(t))-\ell(\gamma(0))}{t}
\]
and the Fisher information tensor is the collection of inner products
\[
g_p(v,w) = E_p d_p\ell(v) d_p\ell(w).
\]

Henceforth we regard $P$ as a Riemannian manifold and assume that $P$ is connected.
Given $p \, d\mu, q \, d\mu \in P$, let $\gamma : [0,1] \rightarrow P$ be a smooth variation such that $\gamma(0)=p$ and $\gamma(1)=q$.  The distance traversed by $\gamma$ is
\[
\mbox{length}(\gamma) = 
\int_0^1 \left[
g_{\gamma(t)} \left( \gamma^\prime(t),\gamma^\prime(t) \right) 
\right]^{1/2} \, dt =
\int_0^1 \left\| \gamma^\prime(t) \right\|_{\gamma(t)} \, dt
\]
and the infimum of these lengths over all such variations defines $i(p,q)$, the {\em information distance}\/
between $p \, d\mu$ and $q \, d\mu$ in $P$.

\subsection{Hellinger Distance}
\label{hell}

\citet[Section 6.8]{murray&rice:1993}
remarked that the fact that the inner products $g_p$ vary with $p$
makes it difficult to discern the global structure of the statistical manifold $P$ directly from Fisher information.  
To remedy this difficulty they defined the square root likelihood,
here denoted $s: P \rightarrow \Re_\Omega$,
by $s(p)=s(p \, d\mu) = 2\sqrt{p}$.
Defining the inner product
\[
\langle f,g \rangle_\mu = \int_\Omega f(x)g(x) \, d\mu(x)
\]
and noting that $2d_ps = sd_p\ell$, we discover that
\begin{eqnarray*}
g_p(v,w) & = & E_p d_p\ell(v)d_p\ell(w) \\
 & = & \int_\Omega \left[ d_p\ell(v) \right](x)
 \left[ d_p\ell(w) \right](x) p(x) \, d\mu(x) \\
 & = & \int_\Omega \left[ \frac{s(p)}{2} d_p\ell(v) \right](x)
 \left[ \frac{s(p)}{2} d_p\ell(w) \right](x) \, d\mu(x) \\
 & = & \int_\Omega \left[ d_ps(v) \right](x)
 \left[ d_ps(w) \right](x) \, d\mu(x) \\
 & = & \left\langle d_ps(v),d_ps(w) \right\rangle_\mu .
\end{eqnarray*}
Hence, if $\gamma$ is a variation in $P$ and $\sigma=s(\gamma)$ is the corresponding variation in $s(P)$, then
\[
\mbox{length}(\gamma) = 
\int_0^1 \left\| \gamma^\prime(t) \right\|_{\gamma(t)} \, dt =
\int_0^1 \left\| \sigma^\prime(t) \right\|_{\mu} \, dt =
\mbox{length}(\sigma) .
\]
The quantity
\begin{equation}
h(p,q) = \left\| s(p)-s(q) \right\|_\mu =
\left\| 2\sqrt{p}-2\sqrt{q} \right\|_\mu
\label{eq:hellinger}
\end{equation}
is the Hellinger distance between the densities $p$ and $q$.
Thus, information distances can be computed by working with Hellinger distance rather than Fisher information.

Now let $\gamma = \{ p_t \, d\mu : t \in (-\epsilon,\epsilon) \}$ denote a smooth variation in the statistical manifold $P$ and consider the Taylor expansion
\begin{equation}
h^2 \left( p_t,p_0 \right) =
h^2 \left( p_0,p_0 \right) +
\left. \frac{d}{dt} h^2 \left( p_t,p_0 \right) \right|_{t=0} t +
\frac{1}{2} \left. \frac{d^2}{dt^2} h^2 \left( p_t,p_0 \right) \right|_{t=0} t^2 + o \left( t^2 \right).
\label{eq:Hellinger1}
\end{equation}
Of course $h^2(p_0,p_0)=0$.  Writing
\begin{eqnarray*}
h^2 \left( p_t,p_0 \right) & = & 
\int_\Omega \left[ 2\sqrt{p_t(x)}-2\sqrt{p_0(x)} \right]^2 \, d\mu(x) \\
 & = & 4 \int_\Omega \left[ p_t(x) - 2\sqrt{p_t(x)p_0(x)} + p_0(x) \right] \, d\mu(x) \\
 & = & 8 - 8 \int_\Omega \left[ p_t(x)p_0(x) \right]^{1/2} \, d\mu(x)
\end{eqnarray*}
and assuming standard regularity conditions that permit differentiation under the integral sign, we obtain
\begin{eqnarray*}
\left. \frac{d}{dt} h^2 \left( p_t,p_0 \right) \right|_{t=0} & = &
-8 \int_\Omega \left. \frac{d}{dt} \left[ p_t(x)p_0(x) \right]^{1/2} \right|_{t=0} \, d\mu(x) \\ & = &
-4 \int_\Omega \left[ p_0(x)p_0(x) \right]^{-1/2} p_0(x)
\left. \frac{d}{dt} p_t(x) \right|_{t=0} \, d\mu(x) \\ & = &
-4 \left. \frac{d}{dt} \int_\Omega p_t(x) \, d\mu(x) \right|_{t=0} \\ 
& = & -4 \left. \frac{d}{dt} 1 \right|_{t=0} =0.
\end{eqnarray*}
Finally,
\begin{eqnarray*}
\left. \frac{d^2}{dt^2} h^2 \left( p_t,p_0 \right) \right|_{t=0} & = &
-8 \int_\Omega \left. \frac{d^2}{dt^2} \left[ p_t(x)p_0(x) \right]^{1/2} \right|_{t=0} \, d\mu(x) \\ & = &
-8 \int_\Omega \left. \frac{d}{dt} \left\{ \frac{1}{2} \left[ p_t(x)p_0(x) \right]^{-1/2} p_0(x) \frac{d}{dt} p_t(x) \right\} \right|_{t=0} \, d\mu(x) \\ & = &
-8 \int_\Omega \left\{ -\frac{1}{4} \left[ p_t(x)p_0(x) \right]^{-3/2} p_0(x) \frac{d}{dt} p_t(x) p_0(x) \frac{d}{dt} p_t(x) + \right. \\ & &
\left. \left. \frac{1}{2} \left[ p_t(x)p_0(x) \right]^{-1/2} p_0(x) \frac{d^2}{dt^2} p_t(x) \right\} \right|_{t=0} \, d\mu(x) \\ & = &
2 \int_\Omega \left[ \frac{\left. \frac{d}{dt} p_t(x) \right|_{t=0}}{p_0(x)} \right]^2 p_0(x) \, d\mu(x) -
4 \int_\Omega \left. \frac{d^2}{dt^2} p_t(x) \right|_{t=0} \, d\mu(x) \\
 & = & 2 \int_\Omega \left. \left[ \frac{d}{dt} \log p_t(x) \right|_{t=0} \right]^2 p_0(x) \, d\mu(x) -  
4 \left. \frac{d^2}{dt^2} \int_\Omega p_t(x) \, d\mu(x) \right|_{t=0} \\
 & = & 2I_\gamma \left(p_0 \right),
\end{eqnarray*}
where $I_\gamma$ denotes Fisher information with respect to the $1$-dimensional submanifold $\gamma$.  Substituting the preceding expressions into (\ref{eq:Hellinger1}) yields
\begin{equation}
h^2 \left( p_t,p_0 \right) =
I_\gamma \left(p_0 \right) t^2 + o \left( t^2 \right).
\label{eq:Hellinger(t)}
\end{equation}
Passing from variations to the (parametrized) manifold $P$,
we write $p_t = p(\cdot,\theta_t)$ and obtain
\begin{equation}
h^2 \left( p \left( \cdot, \theta_t \right), p \left( \cdot, \theta_0 \right) \right) =
\left( \theta_t-\theta_0 \right)^\top I \left( \theta_0 \right) 
\left( \theta_t-\theta_0 \right) + o \left( \left\| \theta_t-\theta_0 \right\|^2 \right).
\label{eq:Hellinger}
\end{equation}
Having derived this expression, the variation $\gamma$ is vestigial and we replace $\theta_t$ in (\ref{eq:Hellinger}) with $\theta$.

\subsection{Minimum Hellinger Distance Estimation}
\label{MHDE}

Following \citep{basu&etal:2011} (with minor changes in notation),
suppose that $x_1,\ldots,x_n \sim p \, d\mu = p(\cdot,\theta) \, d\mu$ and let $\bar{\theta}$ denote the true value of $\theta$.
Let $u(x_i,\theta) = \nabla_\theta \log p(x_i,\theta)$ denote the score function for $P$ and let
\[
Z_n(\theta) = \sqrt{n} \frac{1}{n} \sum_{i=1}^n u \left( x_i,\theta \right).
\]
Under standard regularity conditions, the maximum likelihood estimator $\tilde{\theta}_n$ of $\theta$ is first-order efficient; in particular,
\begin{equation}
\label{eq:MLE1}
\sqrt{n} \left( \tilde{\theta}_n-\bar{\theta} \right) =
I^{-1} \left( \bar{\theta} \right) Z_n \left( \bar{\theta} \right) + o_p(1).
\end{equation}

Let $\hat{p}_n$ denote a nonparametric density estimate of $p$ and define the minimum Hellinger distance estimate (MHDE) of $\theta$ by
\[
\hat{\theta}_n = \arg\min_{\theta \in \Theta}
h \left( p(\cdot,\theta),\hat{p}_n \right) =
\arg\min_{\theta \in \Theta} \int_\Omega
\left[ \sqrt{p(x,\theta)}-\sqrt{\hat{p}_n(x)} \right]^2 \,
d\mu(x).
\]
Under suitable regularity conditions
(see \citep[Section 3.2.2]{basu&etal:2011} and \citep{beran:1977}),
\begin{equation}
\label{eq:MHDE1}
\sqrt{n} \left( \hat{\theta}_n-\bar{\theta} \right) =
I^{-1} \left( \bar{\theta} \right) Z_n \left( \bar{\theta} \right) + o_p(1).
\end{equation}
Thus, both $\tilde{\theta}_n$ and $\hat{\theta}_n$ are first-order efficient estimators.  Typically, $\tilde{\theta}_n$ is more readily  computed and $\hat{\theta}_n$ has better robustness properties.

\section{Information Tests}
\label{IT}

Suppose that $x_1,\ldots,x_n \sim p \, d\mu$, where $p \, d\mu$ lies in the connected $k$-dimensional statistical manifold $P$.  We write $p=p(\cdot,\theta)$, $\bar{p} = p(\cdot,\bar{\theta})$ and test the simple null hypothesis $H_0: p=\bar{p}$ against the composite alternative hypothesis $H_1: p \neq \bar{p}$.
Equivalently, we test $H_0: \theta=\bar{\theta}$ against $H_1: \theta \neq \bar{\theta}$.

Let $\hat{\theta}_n$ denote the MHDE of $\theta$ and
consider the test statistic 
\[
\mbox{ID}_n = i^2 \left( p \left( \cdot,\hat{\theta}_n \right),
p \left( \cdot,\bar{\theta} \right) \right) =
i^2 \left( p \left( \cdot,\hat{\theta}_n \right),
\bar{p} \right),
\]
the squared information distance between $p(\cdot,\hat{\theta}_n)$ and $\bar{p}$ on the statistical manifold $P$.  Because information distance on $P$ behaves locally like Hellinger distance, we begin by studying the local behavior of the related test statistic
\[
\mbox{HD}_n = h^2 \left( p \left( \cdot,\hat{\theta}_n \right),
\bar{p} \right).
\]
Notice that $n \mbox{HD}_n$ differs from the standard Hellinger disparity difference statistic described in
\citep[Section 5.1]{basu&etal:2011},
although it turns out that they are locally asymptotically equivalent.  More precisely, the relation of tests based on HD to Wald tests is analogous to the relation of disparity difference tests to likelihood ratio tests.

We require a technical result about the remainder term in (\ref{eq:Hellinger}).

\begin{lemma}
\label{lm:intersect}
$P(A) \geq 1-\alpha/2$ and $P(B) \geq 1-\alpha/2$ entails $P(A \cap B) \geq 1-\alpha$.
\end{lemma}

\subparagraph{Proof}
Subtracting
\[
P \left( A \cap B \right) +
P \left( A \cap B^c \right) +
P \left( A^c \cap B \right) +
P \left( A^c \cap B^c \right) = 1
\]
from
\begin{eqnarray*}
P \left( A \cap B^c \right) +
P \left( A \cap B \right) +
P \left( A \cap B \right) +
P \left( A^c \cap B \right) & = &
P(A)+P(B) \\ & \geq &
1 - \frac{\alpha}{2} + 1 - \frac{\alpha}{2} =
2 - \alpha
\end{eqnarray*}
yields
\[
P \left( A \cap B \right) \geq
P \left( A \cap B \right) -
P \left( A^c \cap B^c \right) \geq
1 - \alpha.
\]
\hfill $\Box$

\begin{lemma}
\label{lm:remainder}
Let
\[
r(\theta) = h^2 \left( p(\cdot,\theta),\bar{p} \right) -
\left( \theta-\bar{\theta} \right)^\top I \left( \bar{\theta} \right) 
\left( \theta-\bar{\theta} \right).
\]
If (\ref{eq:MHDE1}) holds with $\theta_t=\theta$ and $\theta_0=\bar{\theta}$, then
\[
n \left| r \left( \hat{\theta}_n \right) \right| = o_p(1).
\]
\end{lemma}

\subparagraph{Proof}
Given $c,\alpha>0$, we seek to demonstrate the existence of $N$ such that $n \geq N$ entails
\[
P \left( n \left| r \left( \hat{\theta}_n \right) \right| \geq c \right) < \alpha.
\]
Let $T$ denote the random variable to which $\sqrt{n} (\hat{\theta}_n -\bar{\theta})$ converges in distribution and choose $\epsilon>0$ such that
\[
P \left( \epsilon \| T \|^2 \geq c \right) < \frac{\alpha}{4}.
\]
Choose $N_1$ such that $n \geq N_1$ entails
\[
\left| P \left( \epsilon \left\| \sqrt{n} \left( \hat{\theta}_n -\bar{\theta} \right) \right\|^2 \geq c \right) -
P \left( \epsilon \| T \|^2 \geq c \right) \right| < \frac{\alpha}{4},
\]
and hence that
\[
P \left( B_n^c \right) =
P \left( \epsilon \left\| \sqrt{n} \left( \hat{\theta}_n -\bar{\theta} \right) \right\|^2 \geq c \right) < \frac{\alpha}{4} + \frac{\alpha}{4} =
\frac{\alpha}{2}.
\]
Because $r(\theta) = o(\| \theta-\bar{\theta} \|^2)$,
there exists $\delta>0$ such that $\| \hat{\theta}_n-\bar{\theta} \| < \delta$ entails
\begin{eqnarray*}
\frac{\left| r \left( \hat{\theta}_n \right) \right|}{\left\| \hat{\theta}_n-\bar{\theta} \right\|^2} < \epsilon, & \mbox{ hence } &
n \left| r \left( \hat{\theta}_n \right) \right| < \epsilon \left\| \sqrt{n} \left( \hat{\theta}_n -\bar{\theta} \right) \right\|^2 .
\end{eqnarray*}
Choose $N_2$ such that $n \geq N_2$ entails
\[
P \left( \left\| \hat{\theta}_n-\bar{\theta} \right\| < \delta \right) \geq 1-\frac{\alpha}{2},
\]
hence
\[
P \left( A_n \right) =
P \left( n \left| r \left( \hat{\theta}_n \right) \right|
 < \epsilon \left\| \sqrt{n} \left( \hat{\theta}_n -\bar{\theta} \right) \right\|^2 \right) \geq 1-\frac{\alpha}{2} .
\]
Let $N = \max(N_1,N_2)$.  Then $n \geq N$ entails
\[
P \left( n \left| r \left(\hat{\theta}_n \right) \right|
  < c \right) \geq
P \left( A_n \cap B_n \right) \geq 1-\alpha
\]
by Lemma~\ref{lm:intersect}. \hfill $\Box$

\bigskip

The relation between the HD and Wald statistics is now straightforward.
\begin{theorem}
\label{th:Wald}
Let
\[
W_n = n \left( \tilde{\theta}_n-\bar{\theta} \right)^\top I \left( \bar{\theta} \right) \left( \tilde{\theta}_n-\bar{\theta} \right)
\]
denote the Wald statistic for testing $H_0: \theta=\bar{\theta}$ versus $H_1: \theta \neq \bar{\theta}$.  If (\ref{eq:MLE1}) and (\ref{eq:MHDE1}) hold, then
\[
n\mbox{HD}_n-W_n = o_p(1) .
\] 
\end{theorem}

\subparagraph{Proof}
Applying (\ref{eq:MLE1}), (\ref{eq:MHDE1}), and Lemma \ref{lm:remainder},
\begin{eqnarray*}
n\mbox{HD}_n-W_n & = &
n \left( \hat{\theta}_n-\bar{\theta} \right)^\top I \left( \bar{\theta} \right) \left( \hat{\theta}_n-\bar{\theta} \right) + o_p(1) -
n \left( \tilde{\theta}_n-\bar{\theta} \right)^\top I \left( \bar{\theta} \right) \left( \tilde{\theta}_n-\bar{\theta} \right) \\
 & = & \left[ I^{-1} \left( \bar{\theta} \right) Z_n \left( \bar{\theta} \right) + o_p(1) \right]^\top I \left( \bar{\theta} \right) \left[ I^{-1} \left( \bar{\theta} \right) Z_n \left( \bar{\theta} \right) + o_p(1) \right] + o_p(1) \\
 & & - \left[ I^{-1} \left( \bar{\theta} \right) Z_n \left( \bar{\theta} \right) + o_p(1) \right]^\top I \left( \bar{\theta} \right) \left[ I^{-1} \left( \bar{\theta} \right) Z_n \left( \bar{\theta} \right) + o_p(1) \right] \\
 & = & o_p(1) .
\end{eqnarray*}
\hfill $\Box$

\bigskip

Our Theorem \ref{th:Wald} is analogous to Theorem~1 in \citep{simpson:1989},
which relates a Hellinger deviance test statistic to the likelihood ratio test statistic 
\[
G^2_n = 2 \sum_{i=1}^n \log p \left( x_i, \tilde{\theta}_n \right)/
p \left( x_i, \bar{\theta} \right) .
\]
The asymptotic null distribution of $G^2_n$ and $W_n$ is $\chi^2(k)$;
it follows that the asymptotic null distribution of Simpson's test statistic and our $n\mbox{HD}_n$ is also $\chi^2(k)$.
Furthermore, a contiguity argument (see \citep{simpson:1989}) for details) establishes that these tests have the same asymptotic power at local alternatives of the form $\bar{\theta}+\eta/\sqrt{n}$.
In this sense, our HD test, the Wald test, the likelihood ratio test, and Simpson's Hellinger deviance test are all locally equivalent.

To extend the equivalence to our ID test, we demonstrate that $i^2(p_t,p_0)$ behaves locally like (\ref{eq:Hellinger(t)}).
Recall that a geodesic arc is a variation with zero curvature, hence with constant velocity.
Given $p_0 \in P$, use Lemma 10.3 in \citep{milnor:1963} to choose a neighborhood $W$ of $p_0$ and $\bar{\epsilon}>0$ such that $q \in W$ implies the existence of a unique geodesic variation $\gamma$ connecting $p_0$ and $q$ with $\epsilon=\mbox{length}(\gamma)<\bar{\epsilon}$. It then follows from Theorem 10.4 in \citep{milnor:1963} that $i(q,p_0)=\epsilon$, i.e., that $\gamma$ is the unique path of shortest distance from $p_0$ to $q$.
Parametrizing $\gamma$ by arc length and letting $q=p_\epsilon$, we obtain
\[
\epsilon = i \left( p_\epsilon,p_0 \right) = 
\mbox{length}(\gamma) = 
\int_0^\epsilon \left\| \gamma^\prime(t) \right\|_{\gamma(t)} \, dt
\]
with constant unit velocity
\[
1 = \left\| \gamma^\prime(t) \right\|_{\gamma(t)} =
I_\gamma \left( p_t \right) .
\]
It follows from (\ref{eq:Hellinger(t)}) that
\begin{equation}
h^2 \left( p_\epsilon,p_0 \right) =
I_\gamma \left( p_0 \right) \epsilon^2 + o \left( \epsilon^2 \right) =
\epsilon^2 + o \left( \epsilon^2 \right) =
i^2 \left( p_\epsilon,p_0 \right) + o \left( \epsilon^2 \right) .
\label{eq:d=h}
\end{equation}
Set $\theta_0=\bar{\theta}$.
By arguments analogous to those used to establish Theorem \ref{th:Wald}, we then obtain the following relation.
\begin{theorem}
\label{th:MD=HD}
If (\ref{eq:MHDE1}) holds, then $n\mbox{HD}_n-n\mbox{ID}_n = o_p(1)$.
\end{theorem}
Thus, $\mbox{ID}_n$ and $\mbox{HD}_n$ are locally asymptotically equivalent for testing $H_0: \theta=\bar{\theta}$ versus $H_1: \theta \neq \bar{\theta}$.

Although the information distance, Hellinger distance, Wald, likelihood ratio, and Hellinger disparity distance tests are all locally asymptotically equivalent, only the information distance test attempts to exploit the Riemannian geometry of $P$ when testing nonlocal alternatives.

\section{Restricted Information Tests}
\label{RIT}

Let
$Q = \left\{ p(\cdot,\theta) \, d\mu \; : \; \theta \in \Psi \subset \Theta \right\}$
denote a parametric subfamily of probability distributions in $P$.
Suppose that $Q$
is a $d$-dimensional embedded submanifold of $P$;  equivalently,
suppose that $\Psi$
is a $d$-dimensional embedded submanifold of $\Theta$.
Suppose that $\bar{\theta} \in \Psi$ and that we want to test $H_0 : \theta = \bar{\theta}$, restricting attention to alternatives that lie in $\Psi$.
We emphasize that we are restricting inference to the submanifold, {\em not}\/ testing the null hypothesis that $\theta$ lies in the submanifold.
Two information tests are then available: the unrestricted information test computes information distance on the statistical manifold $P$, whereas the restricted information test computes information distance on the statistical submanifold $Q$.
It is tempting to speculate that restricted information tests are more powerful than unrestricted information tests.

An analogous investigation of restricted likelihood ratio tests was undertaken by \citet{mwt:lrt}, who indeed established that, if 
$d = \mbox{dim}(\Psi) < \mbox{dim}(\Theta) = k$, then the restricted likelihood ratio test is asymptotically more powerful than the unrestricted likelihood ratio test at local alternatives.  As information tests are locally asymptotically equivalent to likelihood ratio tests, they must enjoy the same property.
However, \citet{mwt:lrt} also constructed examples in which the restricted likelihood ratio test is less powerful than the unrestricted likelihood ratio test for certain nonlocal alternatives.  
Unlike restricted likelihood ratio tests,
restricted information tests potentially exploit the {\em global}\/ structure of the statistical submanifold.
This observation motivates investigating the behavior of information tests at nonlocal alternatives.

In what follows we specialize to the case of multinomial distributions, which are widely used (as in \citep{kass:1989}) to illustrate the ideas of information geometry.
Accordingly, consider an experiment with $k+1$ possible outcomes.
The probability model $P=\mbox{Multinomial}(\theta)$ specifies that the outcomes occur with probabilities 
$\theta=(\theta_1,\ldots,\theta_{k+1})$.
It is parametrized by the $k$-dimensional unit simplex in $\Re^{k+1}$,
\[
\Theta = \{ \theta \in [0,1]^{k+1} : \theta_1+\cdots+\theta_{k+1}=1 \},
\] 
or (upon setting $\sigma=\sqrt{\theta}$, defined by setting each $\sigma_i = \sqrt{\theta_i}$)
by that portion of the $k$-dimensional unit sphere that lies in the nonnegative orthant of $\Re^{k+1}$,
\[
\Sigma = \{ \sigma \in [0,1]^{k+1} : \sigma_1^2+\cdots+\sigma_{k+1}^2=1 \}.
\]

One advantage of studying multinomial distributions is the availability of explicit formulas.
If $p=p(\cdot,\theta=\sigma^2)$ and
$q=p(\cdot,\pi=\rho^2)$, then 
\[
h^2(p,q) = \sum_{i=1}^{k+1} \left( 
2\sqrt{\theta_i}-2\sqrt{\pi_i} \right)^2 =
4 \sum_{i=1}^{k+1} \left( \sigma_i-\rho_i \right)^2 =
4 \left\| \sigma-\rho \right\|^2
\]
and we see that Hellinger distance between multinomial distributions corresponds to chordal (Euclidean) distance on $\Sigma$.  Hence, by the law of cosines,
\[
h^2(p,q) = 4 (2 - 2\cos \delta) = 
8-8 \langle \sigma, \rho \rangle,
\]
where $\delta$ is the angle between $\sigma$ and $\rho$.
But $\delta$ is also the great circle (geodesic) distance between $\sigma$ and $\rho$; hence,
\[
i(p,q) = 
2\delta = 2 \arccos \langle \sigma, \tau \rangle,
\]
where the factor of $2$ accrues from (\ref{eq:hellinger}).
It follows that
\[
h^2(p,q) = 8-8 \cos \left( i(p,q)/2 \right),
\]
establishing that the information and Hellinger distances between multinomial distributions are monotonically related.

A second advantage of studying multinomial distributions is that empirical distributions from multinomial experiments are themselves multinomial distributions.  
Suppose that one draws $n$ independent and identically distributed observations from $\mbox{Multinomial}(\theta)$ and counts $\vec{x}=(x_1,\ldots,x_{k+1})$, where $x_i$ records the number of occurrences of outcome $i$.
The empirical distribution of $\vec{x}$ is $\hat{p}_n(\vec{x}) = \vec{x}/n$
and furthermore, because $\vec{x}/n \in \Theta$,
the unrestricted MHDE of $\theta \in \Theta$ is $\hat{\theta}_n(\vec{x}) = \vec{x}/n$.
The restricted MHDE of $\theta \in \Psi$ is
\[
\check{\theta}_n \left( \vec{x} \right) = \arg \min_{\theta \in \Psi}
h^2 \left( \theta,\vec{x}/n \right) = 
\check{\sigma}^2_n \left( \vec{x} \right),
\]
where
\[
\check{\sigma}_n \left( \vec{x} \right) = \arg \max_{\sigma^2 \in \Psi}
\left\langle \sigma, \sqrt{\vec{x}/n} \right\rangle .
\]
Depending on the submanifold $\Psi$,
the calculation of $\check{\theta}_n(\vec{x})$ may require numerical optimization.

Let $i(\cdot,\cdot;\Theta)$ denote information distance on the unrestricted model and let $i(\cdot,\cdot;\Psi)$ denote information distance on the restricted model.
The nonrandomized unrestricted information test with critical value $c_2$ rejects $H_0:\theta=\bar{\theta}$ if and only if
\[
i_n \left( \vec{x};\Theta \right) =
i \left( p \left( \cdot,\hat{\theta}_n \right),
p \left( \cdot,\bar{\theta} \right) ; \Theta \right) =
2 \arccos \left\langle \sqrt{\vec{x}/n}, \sqrt{\bar{\theta}} \right\rangle > c_2.
\] 
The nonrandomized restricted information test with critical value $c_1$ rejects $H_0:\theta=\bar{\theta}$ if and only if
\[
i_n \left( \vec{x};\Psi \right) =
i \left( p \left( \cdot,\check{\theta}_n \right),
p \left( \cdot,\bar{\theta} \right) ; \Psi \right) > c_1.
\]
Because $\vec{x}$ is discrete, randomization may be needed to attain a specified size.
For $n$ sufficiently large, we can use the $1-\alpha$ quantiles $q_{1-\alpha}(k)$ and $q_{1-\alpha}(d)$ of chi-squared distributions with $k$ and $d$ degrees of freedom to select the critical values:
\begin{eqnarray*}
c_2 = \left( q_{1-\alpha}(k)/n \right)^{1/2} & \mbox{ and } &
c_1 = \left( q_{1-\alpha}(d)/n \right)^{1/2}
\end{eqnarray*}
The power functions of the above tests are
\[
\beta_2(\theta) =
P_{\theta \in \Psi} \left( i_n \left( \vec{x};\Theta \right) > c_2 \right)
\]
for the unrestricted information test and
\[
\beta_1(\theta) =
P_{\theta \in \Psi} \left( i_n \left( \vec{x};\Psi \right) > c_1 \right)
\]
for the restricted information test.

\section{Two Trinomial Examples}
\label{HW}

The probability model $\mbox{Trinomial}(\theta)$ specifies that $k+1=3$ outcomes occur with probabilities $\theta=(\theta_1,\theta_2,\theta_3)$.
Define $\psi : [0,1] \rightarrow \Theta$ by 
$\psi(\tau) = \left( \tau^2, 2\tau(1-\tau), (1-\tau)^2 \right)$.
The Hardy-Weinberg subfamily of trinomial distributions is 
parametrized by the embedded submanifold $\Psi = \{ \psi(\tau) : \tau \in [0,1] \}$.
Notice that $\mbox{dim } \Psi =1 < 2 = \mbox{dim } \Theta$.
We write $\mbox{HW}(\tau) = \mbox{Trinomial}(\psi(\tau))$.

Fix $\bar{\tau} \in (0,1)$ and set $\bar{\theta}=\psi(\bar{\tau})$.
We test the simple null hypothesis $H_0:\theta=\bar{\theta}$ against alternatives of the form $\theta=\psi(\tau)$.
The unrestricted information test statistic is
\[
i_n \left( \vec{x};\Theta \right) =
2 \arccos \left( \bar{\tau} \left( x_1/n \right)^{1/2} +
\left[ 2\bar{\tau}(1-\bar{\tau})x_2/n \right]^{1/2} +
(1-\bar{\tau})\left( x_3/n \right)^{1/2} \right).
\]
The restricted MHDE of $\theta \in \Psi$ is
\[
\check{\theta}_n \left( \vec{x} \right) = 
\psi \left( \check{\tau} \left( \vec{x} \right) \right),
\]
where
\[
\check{\tau} \left( \vec{x} \right) = \arg \max_{\tau \in [0.1]}
\left( \tau \left( x_1/n \right)^{1/2} +
\left[ 2\tau(1-\tau)x_2/n \right]^{1/2} +
(1-\tau)\left( x_3/n \right)^{1/2} \right).
\]
Letting $\sigma(\tau) = 2 \psi(\tau)^{1/2}$,
the restricted information test statistic,
$i_n \left( \vec{x};\Psi \right)$,
is computed by integrating
\[
\left\| \sigma^\prime(\tau) \right\| =
2 \left[ 1^2 + \frac{(1-2\tau)^2}{2\tau(1-\tau)} + 1^2 \right]^{1/2}
\]
as $\tau$ varies between $\bar{\tau}$ and 
$\check{\tau}(\vec{x})$.

\begin{table}[tb]
\[
\begin{array}{|cllll|} \hline
x_1,x_2,x_3 & p(\vec{x},\psi(0.3)) &   
i_3 \left( \vec{x};\Theta \right) & \check{\tau}(\vec{x}) &
i_3 \left( \vec{x};\Psi \right) \\ \hline
3,0,0 & 0.000729 &  2.532207  & 1         & 2.803414 \\
2,1,0 & 0.010206 &  1.806363  & 0.8535517 & 1.692687 \\
2,0,1 & 0.011907 &  1.728807  & 1         & 2.803414 \\
1,2,0 & 0.047628 &  1.584191  & 0.7236016 & 1.237653 \\
1,1,1 & 0.111132 &  0.625338  & 0.5       & 0.581973 \\
1,0,2 & 0.064827 &  1.461264  & 0         & 1.639469 \\
0,3,0 & 0.074088 &  1.731487  & 0.5       & 0.581973 \\
0,2,1 & 0.259308 &  0.734627  & 0.2763984 & 0.073708 \\
0,1,2 & 0.302526 &  0.662028  & 0.1464483 & 0.528741 \\
0,0,3 & 0.117649 &  1.590798  & 0         & 1.639469 \\ \hline
\end{array}
\]
\caption{Unrestricted (Trinomial) and restricted (Hardy-Weinberg) information tests of $H_0:\theta=\psi(0.3)$ with $n=3$ observations.  Columns 1--2 list the possible outcomes and their exact probabilities under $H_0$; Column 3 lists the unrestricted information distance of the empirical distributions from the null distribution; Columns 4--5 list the minimum Hellinger distance estimates of the Hardy-Weinberg parameter, $\tau$, and the restricted information distance of the corresponding distributions from the null distribution.}
\label{tbl:HW}
\end{table}

\begin{figure}[tb]
\begin{center}
\includegraphics[width=\textwidth]{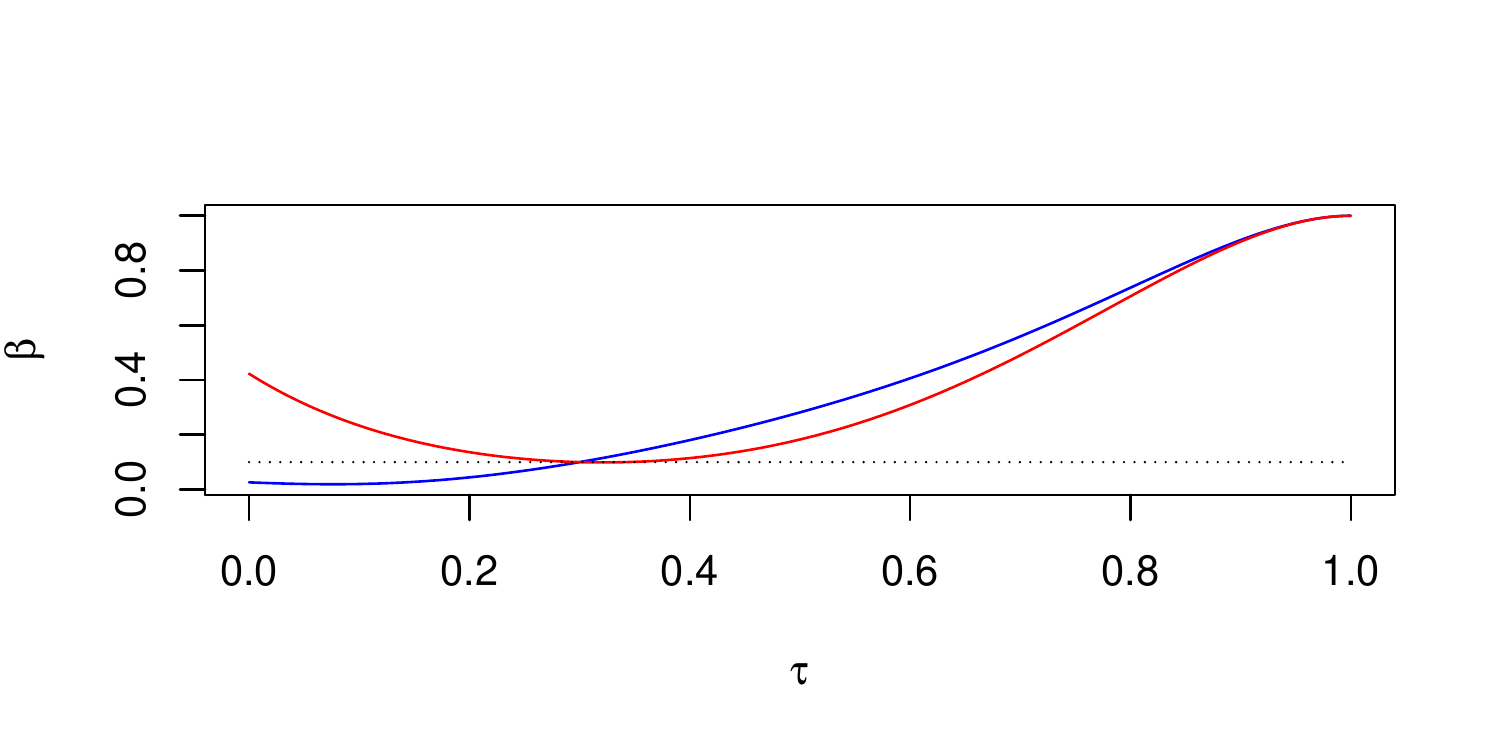} 
\end{center}
\caption{Power of the exact unrestricted ($\beta_2$, plotted in blue) and restricted ($\beta_1$, plotted in red) information tests for testing $H_0: \theta = \psi(0.3)$ with $\alpha=0.1$ (dotted line) and $n=3$.  The alternatives $\{ \theta=\psi(\tau) : \tau \in [0,1] \}$ are displayed on the horizontal axis.  The restricted test is greatly superior for $\tau<0.3$, slightly inferior for $\tau>0.3$.} 
\label{fig:it-hw-3}
\end{figure}

\subparagraph{Example 1}
The trinomial experiment with $n=3$ 
has $10$ possible outcomes, enumerated in the first column of Table~\ref{tbl:HW}.
Consider the unrestricted and restricted information tests of
$H_0:\theta=\psi(0.3)$ with size $\alpha=0.1$.  The exact unrestricted test rejects $H_0$ with certainty if 
\[
C_{2a} = \{(3,0,0),(2,1,0),(0,3,0),(2,0,1)\}
\]
 is observed, and with probability
$(0.1-0.09693)/0.117649 \doteq 0.02609457$ if $C_{2b}=(0,0,3)$ is observed.
The exact restricted test rejects $H_0$ with certainty if 
\[
C_{1a}=\{(3,0,0),(2,0,1),(2,1,0)\}
\] 
is observed, and with probability
$(0.1-0.022842)/0.182476 \doteq 0.4228392$ if 
\[
C_{1b}=\{(1,0,2),(0,0,3)\}
\] 
is observed.  The respective power functions are plotted in
Figure~\ref{fig:it-hw-3}.  The restricted test is dramatically more powerful for $\tau<0.3$, slightly less powerful for $\tau>0.3$.  \hfill $\Box$

\bigskip

The small sample size in Example 1 allows us to illustrate the construction of the unrestricted and restricted information tests, but understates the superiority of the restricted test.  It is curious that the restricted test is {\em less}\/ powerful than the unrestricted test for alternatives $\tau>0.3$, but \citet{mwt:lrt} demonstrated the same anomaly for likelihood ratio tests.
For larger sample sizes, the superiority of the restricted test is unambiguous.

\begin{figure}[tb]
\begin{center}
\includegraphics[width=\textwidth]{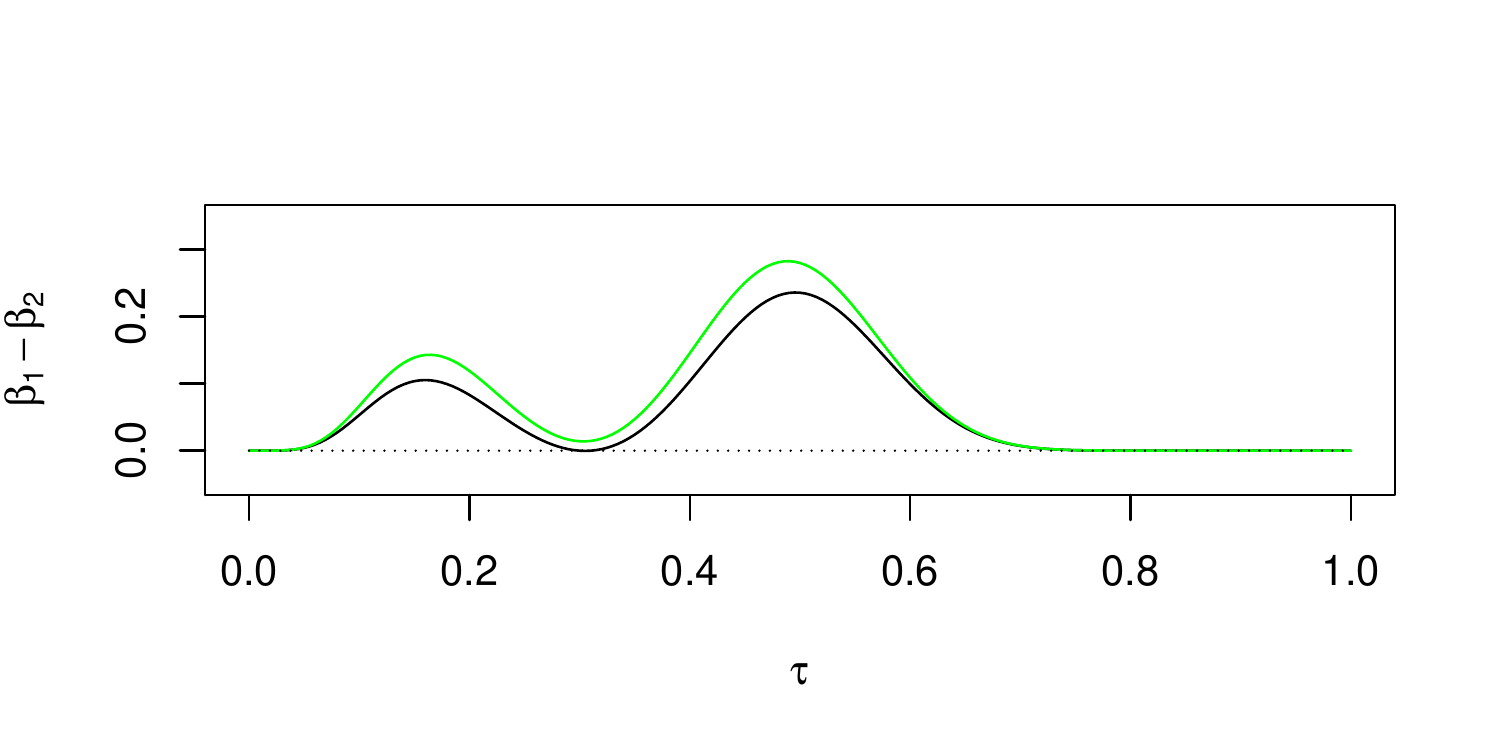} 
\end{center}
\caption{Powers of two restricted information tests ($\beta_1$) minus power of the exact unrestricted information test ($\beta_2$) for testing $H_0: \theta = \psi(0.3)$ with $\alpha=0.05$ and $n=20$.  The black curve corresponds to the exact restricted information test, which has size $\alpha$.  The green curve corresponds to the restricted information test with critical value determined by $\chi^2(1)$, which has size $0.064$. The alternatives $\{ \theta=\psi(\tau) : \tau \in [0,1] \}$ are displayed on the horizontal axis.} 
\label{fig:it-hw-20}
\end{figure}

\subparagraph{Example 2}
The trinomial experiment with $n=20$ has $231$ possible outcomes.
Consider the unrestricted and restricted information tests of
$H_0:\theta=\psi(0.3)$ with size $\alpha=0.05$.  The exact unrestricted test has a critical region of $169$ possible outcomes, with a boundary of one outcome that requires randomization.
The exact restricted test has a critical region of $152$ possible outcomes, with a boundary of one outcome that requires randomization.  The difference in power functions,
$\beta_1(\psi(\tau))-\beta_2(\psi(\tau))$,
is plotted in Figure~\ref{fig:it-hw-20}.
The restricted test is clearly superior, although careful examination reveals that it is slightly inferior for alternatives slightly greater than 0.3.  For example, 
\[
\beta_1(\psi(0.305))-\beta_2(\psi(0.305)) \doteq
-0.0002842388.
\]
For comparison, a $\chi^2(1)$ approximation yields a critical value of $c_1 = 0.4382613$.  The corresponding critical region is slightly larger than the exact critical region, containing an additional $5$ outcomes.  Using a larger critical region increases the probability of rejection, in particular to a size of $0.06402558$.  This power function, minus $\beta_2(\psi(\tau))$, is also plotted in Figure~\ref{fig:it-hw-20}.  \hfill $\Box$

\section{Approximate Information Tests}
\label{AIT}
So far, our exposition has glossed the computational challenges posed by information tests.  For multinomial manifolds, the empirical distributions lie on the manifold and information distance can be computed by a simple formula.  For the $1$-dimensional Hardy-Weinberg submanifold, minimum Hellinger distance estimates require numerical optimization, geodesic variations are apparent by inspection, and computing an information distance requires numerical integration.  In general, however, the information tests described in Sections \ref{IT} and \ref{RIT} necessitate overcoming the following challenges:
\begin{enumerate}

\item  Numerical optimization on the submanifold to determine the minimum Hellinger distance estimate, $\check{\theta}_n$.

\item  Determining the geodesic variation between $\check{\theta}_n$ and the hypothesized $\bar{\theta}$.  If the submanifold is $1$-dimensional, then this is easily accomplished by inspection;
if $d>1$, then the geodesic variation must be determined by solving a potentially intractable problem in the calculus of variations.  

\item  Numerical integration along the geodesic variation to determine the information distance between $\check{\theta}_n$ and $\bar{\theta}$.

\end{enumerate}
We now propose procedures that circumvent these challenges.
The key idea that underlies these procedures is that information distance is locally approximated by Hellinger distance.

In what follows, we assume that the problems described above are difficult or intractable, but that we can identify a finite set of distributions in the submanifold $Q = \left\{ p(\cdot,\theta) \, d\mu \; : \; \theta \in \Psi \subset \Theta \right\}$.  For example, in the case of the Hardy-Weinberg submanifold, we might identify $m$ trinomial distributions by drawing $\tau_1,\ldots,\tau_m \sim \mbox{Uniform}(0,1)$.  Combined with the hypothesized distribution, we thus have $m+1$ distributions in $Q$ from which we hope to learn enough about the Riemannian structure of $Q$ to approximate the methods of Section~\ref{RIT}.

Elaborating on Figure~\ref{fig:AIT},
we propose the following procedure for testing $H_0: \theta = \bar{\theta}$.
\begin{enumerate}

\item Identify $\theta_1,\ldots,\theta_m \in \Psi$ and compute the  $(m+1)m/2$ pairwise Hellinger distances $h_{ij}$ between the $\bar{p},p_1,\ldots,p_m$ that correspond to 
$\bar{\theta},\theta_1,\ldots,\theta_m$.

\item Use the pairwise Hellinger distances to form ${\mathcal G}$, a graph with $m+1$ vertices corresponding to the $m+1$ distributions.  Connect vertices $i$ and $j$ when $h_{ij}$ is sufficiently small, so that ${\mathcal G}$ localizes the structure of the submanifold $Q$.  
Weight edge $i \leftrightarrow j$ by $h_{ij}$.

This is a standard construction in manifold learning, e.g.,
\citep{isomap:2000,lle:2000},
although our application of manifold learning techniques to statistical rather than data manifolds appears to be novel.
The most popular constructions are either (a) connect $i$ and $j$ if and only if $h_{ij} \leq \epsilon$, or (b) connect $i$ and $j$ if and only if $i$ is a $K$-nearest neighbor (KNN) of $j$ or $j$ is a KNN of $i$.  The choice of the localization parameter ($\epsilon$ or $K$) is a model selection problem.
It is imperative that the localization parameter be chosen so that 
${\mathcal G}$ is connected.

\item Compute $\Delta = [ \delta_{ij} ]$, the $(m+1) \times (m+1)$ dissimilarity matrix of pairwise shortest path distances in 
${\mathcal G}$.

Here we appropriate the key idea of the popular manifold learning procedure Isomap \citep{isomap:2000}.  A path in ${\mathcal G}$ is a discrete approximation of a variation in $Q$.
The length of a path is the sum of its Hellinger distance edge weights, hence a discrete approximation of the integral that defines the length of the approximated variation.
The shortest path between vertices $i$ and $j$ approximates the geodesic variation between distributions $i$ and $j$, hence the shortest path distance $\delta_{ij}$ approximates the information
distance between distributions $i$ and $j$.

\item For a suitable choice of $r$, embed $\Delta$ in $\Re^r$ by minimizing a suitably weighted raw stress criterion,
\[
\sigma(Z) = \sum_{i<j} w_{ij} \left[
\left\| z_i-z_j \right\| - \delta_{ij} \right]^2,
\]
where the coordinates of $z_i \in \Re^r$ appear in row $i$ of the $(m+1) \times r$ configuration matrix $Z$.

Isomap \citep{isomap:2000} embeds shortest path distances by classical multidimensional scaling
\citep{torgerson:1952,gower:1966},
which minimizes a squared error criterion for pairwise inner products.  The widely used raw stress criterion is more directly related to our objective of modeling shortest path distance with Euclidean distance; it also provides greater flexibility through its ability to accommodate different weighting schemes.
The raw stress criterion can be numerically optimized by majorization \citep{deleeuw:1988}, several iterations of which usually provides a useful embedding, or by Newton's method 
\citep{kearsley&etal:newton}, which has better local convergence properties.

The choice of $r$ is a model selection problem.  While $r=d$ is nearly universal in conventional manifold learning, $r>d$ may provide a more faithful Euclidean representation of the geodesic structure of $Q$.

\item  From $x_1,\ldots,x_n \sim p$, construct a nonparametric density estimate $\hat{p}_n$.  Compute the Hellinger distances of 
$\hat{p}_n$ from $p_1,\ldots,p_m$ and let
$j_1,\ldots,j_\ell$ index the nearest $\ell \geq r$ distributions.
Embed $\hat{p}_n$ in the previously constructed representation by a suitable out-of-sample embedding technique.  Let $y(\vec{x}) \in \Re^r$ denote the resulting representation of 
$\hat{p}_n$.  The proposed approximate information test rejects $H_0: \theta=\bar{\theta}$ if and only if the test statistic
\[
\hat{i}_n \left( \vec{x};\Psi \right) =
\left\| y \left( \vec{x} \right) -\bar{z} \right\|,
\]
where $\bar{z}$ corresponds to $\bar{\theta}$,
is sufficiently large.

A comprehensive discussion of how to embed $\hat{p}_n$ using only its $\ell$ nearest neighbors is beyond the scope of this manuscript.  For $r=1$ and $\ell=2$, one can use the law of cosines to project $\hat{p}_n$ into the line that contains $z_{j_1}$ and $z_{j_2}$.  This construction is a special case of out-of-sample embedding into a principal components representation.  See 
\cite{gower:1968}
for a general formula that uses pairwise squared distances; see 
\cite{williams&seeger:2001}
for a general formula that uses pairwise inner products.
For the simulations in Example~4, 
we simply set $y(\vec{x})$ equal to the centroid of
$z_{j_1},\ldots,z_{j_\ell}$.

\item Estimate a significance probability by generating simulated random samples $\vec{x}_i$ of size $n$ from the hypothesized distribution $\bar{p}$.
Perform the previous step for each $\vec{x}_i$ and compute the fraction of $\vec{x}_i$ for which
\[
\hat{i}_n \left( \vec{x}_i;\Psi \right) \geq
\hat{i}_n \left( \vec{x};\Psi \right).
\]

\end{enumerate}

\begin{figure}[tb]
\begin{center}
\includegraphics[width=0.9\textwidth]{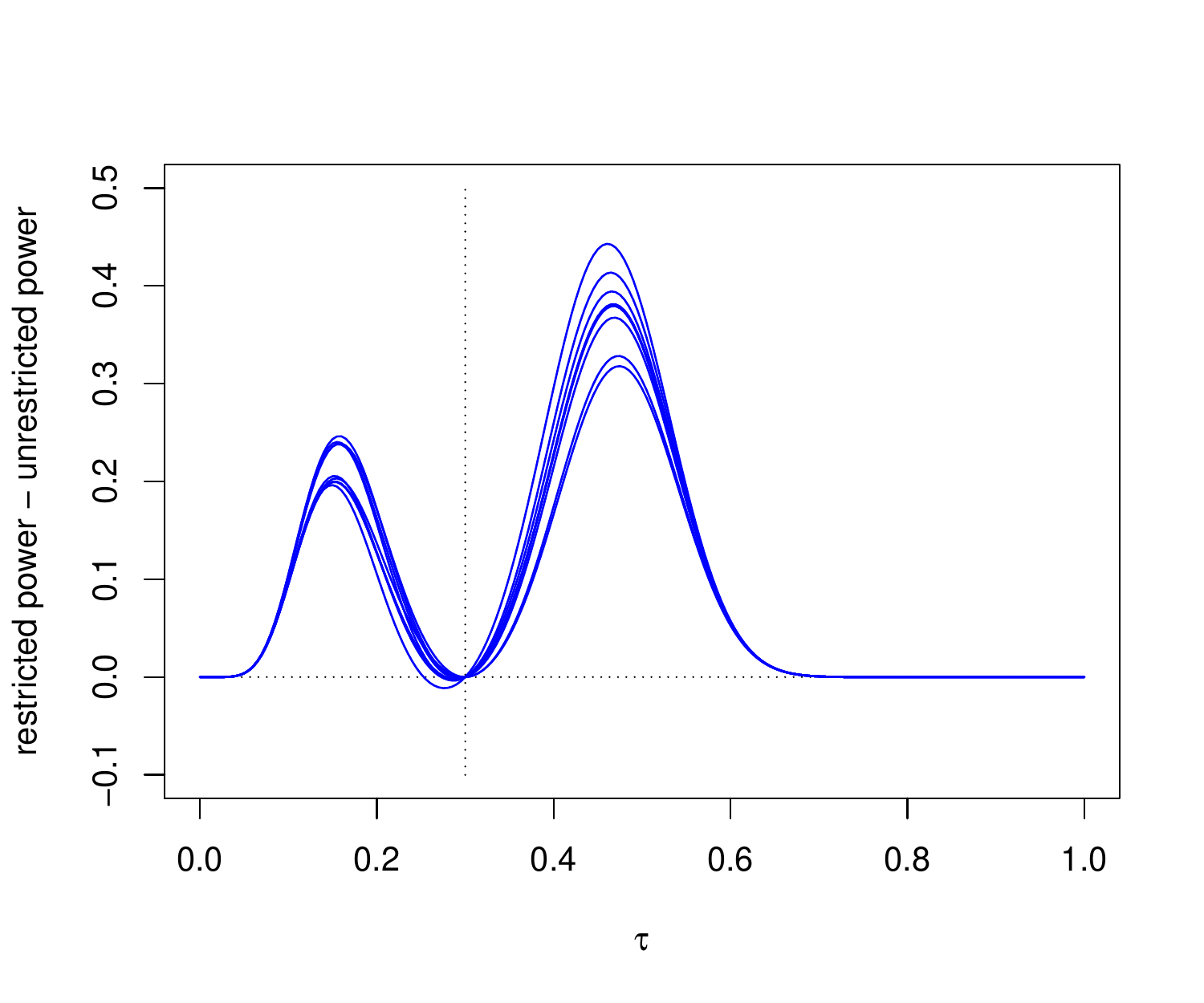} 
\end{center}
\caption{Powers of ten approximate restricted information tests minus power of the exact unrestricted information test for testing $H_0: \theta = \psi(0.3)$ versus $H_1: \theta \in \{ \psi(\tau) : \tau \in [0,1] \}$
with $\alpha=0.05$ and $n=30$.  Each test was randomized to have size $\alpha$.
Each restricted test was constructed using only a random sample of $m=9$ points from the Hardy-Weinberg submanifold.} 
\label{fig:it-hw30-n30-m09-exact}
\end{figure}

\subparagraph{Example 3}
As in Section~\ref{HW}, we consider the Hardy-Weinberg submanifold of  
$\mbox{Trinomial}(\theta)$, defined by $\psi(\tau) = \left( \tau^2, 2\tau(1-\tau), (1-\tau)^2 \right)$ for $\tau \in [0,1]$.  Using $n=30$ trials, we test $H_0: \theta = \psi(0.3)$ by two methods:
\begin{itemize}

\item[a] The information test on the unrestricted manifold of trinomial distributions, for which information distance can be computed by explicit calculation.

\item[b] Ten approximate information tests on estimated $1$-dimensional submanifolds, each constructed using $\bar{\tau}=0.3$ and $\tau_1,\ldots,\tau_9 \sim \mbox{Uniform}[0,1]$.  Shortest path distances on 5NN graphs weighted by pairwise Hellinger distances were embedded in $\Re$ using the unweighted raw stress criterion.  Empirical distributions were then embedded by applying the law of cosines to the $\ell=2$ nearest neighbors.  

\end{itemize}
In each case, a randomized test was constructed to have size $\alpha=0.05$.  
Note that we use the adjectives {\em exact}\/ and {\em approximate}\/ to indicate whether the information distance was computed exactly or approximated by random sampling and manifold learning, not to describe the size of the test.

The power function of the exact unrestricted test was subtracted from the power functions of the ten approximate restricted tests, resulting in the ten difference functions displayed in Figure~\ref{fig:it-hw30-n30-m09-exact}.  Except occasionally for values of $\tau$ slightly less than $0.3$, the approximate restricted tests are consistently more powerful than the exact unrestricted test---often dramatically so.    \hfill $\Box$

\bigskip

We now return to the Motivating Example in Section~\ref{intro} and illustrate the proposed methodology.

\subparagraph{Example 4}
We parametrize the family of multinomial distributions with $7$ possible outcomes by $\Sigma$, the portion of the $6$-dimensional unit sphere in $\Re^7$ that lies in the nonnegative orthant.  The null hypothesis to be tested is 
\[
H_0: \sigma = \bar{\sigma} = (0.3,0.3,0.3,0.5,0.4,0.4,0.4).
\]
Define $\psi : [0,\pi/2]^2 \rightarrow \Sigma$ by
\[ 
\psi(\tau) = \left( 0.3,0.3,0.3,0.5,
\rho \cos \tau_1 \sin \tau_2,
\rho \sin \tau_1 \sin \tau_2,
\rho \cos \tau_2 \right),
\]
where $\rho^2 = 0.48$.
The $2$-dimensional subfamily of multinomial distributions  
defined by the embedded submanifold $\Psi = \{ \psi(\tau) : \tau \in [0,\pi/2]^2 \}$ is a spherical subfamily in the sense of \cite{gous:1999}.
Notice that setting $\tau_1 = \pi/4$ and $\tau_2 = \arctan \sqrt{2}$
results in $\psi(\tau) = \bar{\sigma}$. 

We want to test $H_0$ against alternatives that lie in $\Psi$.
If $\Psi$ was known, then we could perform a restricted likelihood ratio test.
The likelihood of $o=(3,5,4,6,9,2,1)$ under $\sigma = \psi(\tau)$ is
\[
L_o(\psi(\tau)) = C \cdot 0.09^{3+5+4} \cdot 0.25^6 \cdot
\left( \rho \cos \tau_1 \sin \tau_2 \right)^{2 \cdot 9} \cdot
\left( \rho \sin \tau_1 \sin \tau_2 \right)^{2 \cdot 2} \cdot
\left( \rho \cos \tau_2 \right)^{2 \cdot 1}.
\]
To find the restricted maximum likelihood estimate of $\tau$,
it suffices to minimize
\[
f(\tau) = \left( -18 \log \cos \tau_1 - 4 \log \sin \tau_1 \right) +
\left( -2 \log \cos \tau_2 - 22 \log \sin \tau_2 \right) =
f_1 \left( \tau_1 \right) + f_2 \left( \tau_2 \right)
\]
subject to simple bound constraints $\tau \in [0,\pi/2]^2$.
The objective function $f$ is separable: it suffices to choose $\tau_1$ to minimize $f_1$ and $\tau_2$ to minimize $f_2$.
Furthermore,  $f_1$ and $f_2$ are each strictly convex on $[0,\pi/2]$ (each has a strictly positive second derivative on $(0,\pi/2)$), with unique global minimizers at
\begin{eqnarray*}
\breve{\tau}_1 = \arcsin \sqrt{2/11} \doteq 0.4405107 & \mbox{ and } &
\breve{\tau}_2 = \arcsin \sqrt{11/12} \doteq 1.277954.
\end{eqnarray*}
The restricted likelihood ratio test statistic is then
\begin{eqnarray*}
-2 \log L_o \left( \bar{\sigma} \right) / 
L_o \left( \psi \left( \breve{\tau} \right) \right) & = &
-2 \log 0.16^{12} / \left( 0.36^9 \cdot 0.08^2 \cdot 0.04 \right) \\
 & = & 36 \log 3 - 44 \log 2 \\
 & \doteq &  9.051566.
\end{eqnarray*}
The standard asymptotic approximation of the null distribution of the test statistic is a chi-squared distribution with $2$ degrees of freedom, resulting in 
an approximate significance probability of $\mbox{\bf p} = 0.01082623$.
This significance probability is considerably smaller than the significance probabilities that resulted from the unrestricted Pearson and likelihood ratio tests performed in the Motivating Example.  Unlike them, it causes rejection of $H_0$ at significance level $\alpha =0.05$.

\begin{figure}[tb]
\begin{center}
\includegraphics[width=0.7\textwidth]{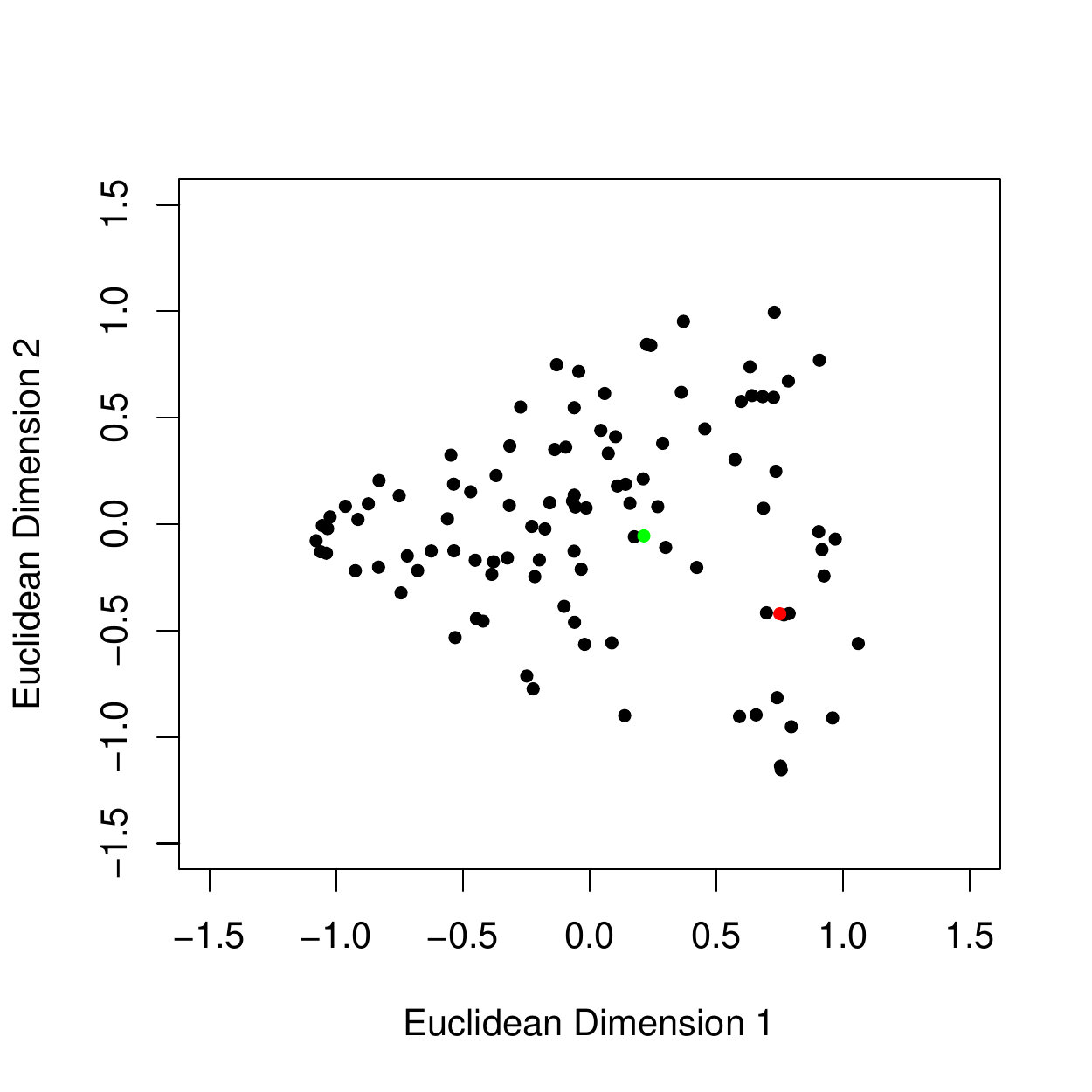} 
\end{center}
\caption{The estimated submanifold in Example~4.
The $m=100$ possible distributions generated by sampling are indicated by
\dot; the null hypothesis is indicated by {\gdot}; and
the minimum distance estimate based on the empirical distribution is indicated by {\rdot}.
The proposed test statistic is $\|${\rdot}$-${\gdot}$\|$, which leads to an
estimated significance probability of $0.0275$.} 
\label{fig:it2d}
\end{figure}

Of course, it is only possible to perform a likelihood ratio test of $H_0 : \sigma = \bar{\sigma}$ versus $H_1: \sigma \in \Psi$ if $\Psi$ is known.  We are concerned with the case that $\Psi$ is unknown, but elements of $\Psi$ can be obtained by sampling.  To simulate that scenario, we drew $\tau_1,\ldots,\tau_{100} \sim \mbox{Uniform}[0,\pi/2]^2$ and computed $\sigma_i = \psi(\tau_i)$.  As reported in Section~\ref{intro}, the first two principal components of the corresponding $\theta_i$ account for $96$\% of the variation in the $m=100$ multinomial parameter values.  The vectors $\bar{\sigma},\sigma_1,\ldots,\sigma_m \in \Re^7$ were then embedded in $\Re^2$ by the manifold learning procedure described above, resulting in Figure~\ref{fig:it2d}.  In this representation of the estimated submanifold, $\hat{\Psi}$, $\sigma_1,\ldots,\sigma_m$ are indicated by \dot, $\bar{\sigma}$ is indicated by {\gdot}, and
$y(\vec{x})$ is indicated by {\rdot}.
Repeating this procedure on $2000$ simulated samples of size $n=30$ drawn from the null distribution resulted in just $55$ larger values of the test statistic, i.e., the estimated significance probability is $55/2000=0.0275$.
The evidence against $H_0$ produced by the restricted approximate information test is slightly less compelling than the evidence produced by the restricted likelihood ratio test (for which $\Psi$ is known), but is more compelling than the unrestricted Pearson or likelihood ratio tests.    \hfill $\Box$

\section{Discussion}
\label{disc}

It is widely believed throughout the statistics community that restricted tests are more powerful than unrestricted tests.  Indeed, although restricted tests may not be uniformly more powerful than unrestricted tests, our experience has been that the former generally outperform the latter.  In consequence, we generally prefer restricted likelihood ratio tests to unrestricted likelihood ratio tests.  But restricted likelihood ratio tests can only be constructed when the restriction to a parametric family of probability distributions is known and tractable.  It is not clear that the low-dimensional structure of a restricted submanifold of distributions can be exploited when the submanifold is unknown.

For simple null hypotheses, we have proposed information tests that are locally asymptotically equivalent to likelihood ratio.  Except in the special case of $1$-dimensional submanifolds, these tests are computationally less tractable than likelihood ratio tests---typically intractable.  Unlike likelihood ratio tests, however, information tests can be approximated when the relevant submanifold of distributions is unknown.  

While local asymptotic theory commends the use of restricted tests, it does not guarantee that finite approximations of restricted tests will outperform unrestricted tests using finite sample sizes.  Nevertheless, we report examples in which the unknown submanifold of distributions can be estimated well enough to realize gains in power.
A preliminary version of our methodology has already been used
to infer brainwide neural-behavioral maps from optogenetic experiments on {\em Drosophila}\/ larvae
\citep{vogelstein&etal:2014}.

A natural extension of the methods reported herein will be to the case obtained in \citep{vogelstein&etal:2014}, in which the randomly generated 
$\theta_1,\ldots,\theta_m \in \Psi$ are replaced by randomly generated 
$\vec{\theta}_1,\ldots,\vec{\theta}_m$ near $\Psi$.  The same methods can be used (and we have used them successfully), but replacing known $\theta_i$ with approximated $\vec{\theta}_i$ introduces another layer of uncertainty.
We are currently exploring such extensions in related work.

\section*{Acknowledgments}
This work was partially supported by DARPA 
XDATA contract FA8750-12-2-0303,
SIMPLEX contract N66001-15-C-4041,
GRAPHS contract N66001-14-1-4028,
and D3M contract FA8750-17-2-0112.

\bibliography{$HOME/lib/tex/stat,$HOME/lib/tex/mds,$HOME/lib/tex/math,$HOME/lib/tex/mwt,$HOME/lib/tex/bio}

\end{document}